\documentclass{ifacconf}
\usepackage[utf8]{inputenc}
\usepackage{amsmath,amssymb,amsfonts}
\usepackage{algcompatible}
\usepackage{algorithm}
\usepackage{graphicx}
\usepackage{textcomp}
\usepackage{booktabs}
\usepackage{subcaption}
\usepackage{multirow}
\usepackage{natbib}
\usepackage{comment}
\usepackage{xcolor}
\usepackage{url}    
\usepackage[T1]{fontenc}
\usepackage{lmodern}

\usepackage[inline]{enumitem}
\usepackage[font={small},labelfont={bf}]{caption}

\usepackage{xcolor}

\usepackage{setspace}

\DeclareMathOperator*{\argmin}{argmin}

\usepackage[
  separate-uncertainty = true,
  multi-part-units = repeat,
  parse-numbers=false
]{siunitx}

\begin{document}

\begin{frontmatter}

\title{Parallel Dynamic Programming for\ Conic Linear Quadratic Control
\thanksref{footnoteinfo}} 

\thanks[footnoteinfo]{This work is partially supported by NWO (Project AMADeuS) and US NSF (Award 2411369). Opinions, findings, conclusions, or recommendations expressed in this material are those of the authors and do not necessarily reflect those of the funding organizations.}

\author[First]{Luyao Zhang} 
\author[Second]{Gabriel Bravo-Palacios} 
\author[Second]{Brian Plancher}
\author[First]{Sergio Grammatico}

\address[First]{Delft Center for Systems and Control, TU Delft, The Netherlands (e-mail: \{l.zhang-7, s.grammatico\}@tudelft.nl).}
\address[Second]{Barnard College, Columbia University, 
   New York, NY 10027 USA \& Dartmouth College, 
   Hanover, NH 03755, USA \\(e-mail: \{gbravo, plancher\}@dartmouth.edu)}

\begin{abstract}                
% Abstract of 50--100 words
Linear Quadratic (LQ) control problems are at the heart of linear control theory and Model Predictive Control (MPC).
While performant, standard approaches to solving such problems are inherently serial, limiting real-time scalability despite the parallel computing power available on modern multi-core CPUs.
Contributing to addressing this challenge and motivated by ``divide and conquer'' strategies, we present a parallel-in-time approach that solves computationally demanding conic optimal control problems through the use of the alternating direction method of multipliers (ADMM). 
In particular, we formulate the inner primal update of ADMM as an LQ problem and split the reformulated problem along the time horizon. This enables us to derive a variant of the Riccati recursion using dynamic programming to solve each subproblem in parallel. 
Numerical benchmarks on two real-world applications demonstrate as much as a 5x speedup compared to existing related approaches on multi-core CPU hardware.
\end{abstract}


\end{frontmatter}
%===============================================================================

%%%%%%%%%%%%%%%%%%%%%%%%%%%%%%%%%%%%%%%%%%%%%%%%%%%%%%%%%%%
%%%%%%%%%%%%%%%%%%%%%%%%%%%%%%%%%%%%%%%%%%%%%%%%%%%%%%%%%%%
%%%%%%%%%%%%%%%%%%%%%%%%%%%%%%%%%%%%%%%%%%%%%%%%%%%%%%%%%%%
%%%%%%%%%%%%%%%%%%%%%%%%%%%%%%%%%%%%%%%%%%%%%%%%%%%%%%%%%%%

\section{Introduction}
Model Predictive Control (MPC) has become increasingly popular in practice due to its ability to adapt to dynamic environments and systematically handle complex constraints, e.g., second-order conic constraints that commonly appear in robotics and aerospace control problems involving friction~\citep{WensingOptimization2024} and thrust limits~\citep{TracyManchester2021_LowThrustKS}. 
These algorithms operate by repeatedly solving finite-horizon optimal control problems online.
Fueled by advances in numerical optimization solvers, real-time conic nonlinear MPC is now practical for many robotics tasks and modalities~\citep{WensingOptimization2024}.

In most efficient MPC implementations, the key computational primitive is the Linear-Quadratic (LQ) problem~\citep{rao_application_1998,Vanroye_FATROP_2023}. This includes both historic and modern iterative schemes such as Sequential Quadratic Programming (SQP)~\citep{betts2010practical,jordana_stagewise_2025,Verschueren2021} and Differential Dynamic Programming (DDP)~\citep{jacobson1970differential,LiUnified2023,HowellALTRO2019}. 
In most of these approaches, the resulting LQ (sub)problem can be interpreted as a quadratic optimization problem subject to linear dynamics constraints. 
Its corresponding KKT system exhibits a banded structure, which can be solved using a sparse $LDL^\top$ factorization, or an efficient Riccati recursion~\citep{dunn_efficient_1989} inspired by dynamic programming. 
We note that, as discussed in Jordana et al.~(\citeyear{jordana_stagewise_2025}), the two approaches are closely related, and many variants, such as the square-root Riccati recursion~\citep{Frison_LQ_2013}, have been proposed for further efficiency.

However, the complexity of the aforementioned methods scales linearly with the length of the prediction horizon and becomes a bottleneck in long-horizon problems~\citep{BettsLowThrust2003}. Moreover, while current solvers for (conic) trajectory optimization are serial~\citep{garstka_cosmo_2021, odonoghue_conic_2016}, recent advances in parallel computing motivate the development of new parallel algorithms for such problems. 
One key class of such parallel algorithms is ``divide and conquer'' approaches, which decompose
the problem into smaller subproblems that can be solved
independently, followed by a consensus step that combines their partial solutions.
In particular, Wright~(\citeyear{wright_partitioned_1991}) proposes two such methods, partitioned dynamic programming (PDP) and partitioned Riccati recursion (PRI), which differ in their treatment of the system of equations associated with the LQ problem, and in the way they parameterize the subproblems. There are multiple variants and extensions to these approaches~\citep{nielsen_parallel_2015,laine_parallelizing_2019,SarkkaTemporal2023}; in particular, Jallet et al.~(\citeyear{jallet_parallel_2024}) consider general LQ problems with implicit dynamics and stage-wise equality constraints, and
under explicit dynamics, apply a reduction phase that coincides with that of the PDP method. 
Table \ref{tab:comparison_parallel_methods} summarizes the comparison among these various methods.

\begin{table*}[h!]
\caption{Comparison of parallel methods.}
\label{tab:comparison_parallel_methods}
\centering
\resizebox{0.9\linewidth}{!}{%
\begin{tabular}{lcccc}
\toprule
\multirow{2}{*}{\textbf{Method}} & 
\multicolumn{2}{c}{\textbf{Reduction}} & 
\multicolumn{2}{c}{\textbf{Consensus}} \\ 
\cmidrule(lr){2-3} \cmidrule(lr){4-5}
 & Factorization Type & Computational Efficiency & Factorization Type & Computational Efficiency \\
\midrule
\texttt{PDP} \citep{wright_partitioned_1991} & Stage $LL^\top$ & Medium & Banded $LDL^\top$ & Medium \\
\texttt{PRI} \citep{wright_partitioned_1991} & Stage $LU$ & Low & Banded $LDL^\top$ & Medium \\
\citealp{nielsen_parallel_2015} & Stage $LL^\top$ & Medium & Stage $LL^\top$ & Low \\
\texttt{ProxLQR} \citep{jallet_parallel_2024} & Stage $LDL^\top$  & Medium & Stage $LDL^\top$ & High \\
\midrule
\texttt{PDPLQR} (ours) & Stage square-root $LL^\top$ & High & Stage $LU/LL^\top$ & High \\
\bottomrule
\end{tabular}%
}
\end{table*}

Despite the advances achieved, with the exception of Jallet et al.~(\citeyear{jallet_parallel_2024}), the aforementioned methods do not account for constraints at all, let alone conic constraints.
We build on this previous work and introduce a new LQ solver which is not only more computationally efficient through its parallel-in-time approach, but can also solve large-scale conic optimal control problems through the alternating direction method of multipliers (ADMM) for constraint handling. 
In particular, we decompose the ADMM primal problem into parallel fixed-end LQ subproblems, which we solve through our own variant of the Riccati recursion, achieving increased computational efficiency. 
Compared to related and common approaches for solving such problems, we demonstrate as much as a 5x speedup on multi-core CPU hardware. We release our software open-source at: \texttt{\url{https://github.com/Luyao787/PDP-LQR}}.

%%%%%%%%%%%%%%%%%%%%%%%%%%%%%%%%%%%%%%%%%%%%%%%%%%%%%%%%%%%
%%%%%%%%%%%%%%%%%%%%%%%%%%%%%%%%%%%%%%%%%%%%%%%%%%%%%%%%%%%
%%%%%%%%%%%%%%%%%%%%%%%%%%%%%%%%%%%%%%%%%%%%%%%%%%%%%%%%%%%
%%%%%%%%%%%%%%%%%%%%%%%%%%%%%%%%%%%%%%%%%%%%%%%%%%%%%%%%%%%

\section{Background} \label{sec:problem_formulation}
\subsection{Dense Linear Algebra}
This work uses the following dense linear algebra operations: matrix-matrix/vector multiplication $(\texttt{gemm}/\texttt{gemv})$, triangular matrix-matrix/vector multiplication $(\texttt{trmm}/\texttt{trmv})$, symmetric rank-$k$ update $(\texttt{syrk})$, Cholesky factorization, triangular solve $(\texttt{trsv})$, as well as LU factorization and solve.
More technical details can be found in \citep[Appendix C]{Boyd_Vandenberghe_2004}.

\subsection{The Conic LQ Problem}
We consider the following convex conic LQ problem:
\begin{subequations} \label{eq:conic_lqr}
\begin{align} 
\min_{\boldsymbol{x}, \boldsymbol{u}, \boldsymbol{s}}\; & \sum_{k=0}^{N-1}\ell_k(x_k,u_k) + \ell_N(x_N) \label{eq:cost} \\[2pt] 
\text{s.t.}\quad & x_{k+1} = A_k x_k + B_k u_k + c_k, \label{eq:lin_dyn}\\
& D_{x,k} x_k + D_{u,k} u_k + s_k = e_k,\\ \label{eq:r_constraint}
& s_k \in \mathcal{K}_k, ~~~~~~\qquad k = [0,N-1],\\
& D_{x,N} x_N + s_N = e_N, \quad s_N \in \mathcal{K}_N, \label{eq:t_constraint}
\end{align}
\end{subequations}
where $N$ is the prediction horizon. The decision variables include the state and control sequences, $\boldsymbol{x} := \left( x_1, \dots, x_{N} \right)$ with $x_k \in \mathbb{R}^{n_x}$, and $\boldsymbol{u} := \left( u_0, \dots, u_{N-1} \right)$ with $u_k \in \mathbb{R}^{n_u}$, respectively, and the slack variables $\boldsymbol{s} := \left( s_0, \dots, s_{N-1} \right)$ with $s_k \in \mathbb{R}^{n_{c,k}}$. 
The stage cost $\ell_k(\cdot)$ and the terminal cost $\ell_N(\cdot)$ in~\eqref{eq:cost} have the quadratic form:
\begin{align*}
    \ell_k(x_k, u_k) &= 
    \begin{bmatrix} q_{k} \\  r_{k}     \end{bmatrix} ^\top     
    \begin{bmatrix} x_k \\ u_k \end{bmatrix} + \frac{1}{2}    
    \begin{bmatrix} x_k \\ u_k \end{bmatrix} ^\top     
    \begin{bmatrix} Q_{k} &  M_{k}^\top\\ M_{k} & R_{k}    \end{bmatrix}     
    \begin{bmatrix} x_k \\ u_k \end{bmatrix}, \\[1pt]
    \ell_N(x_N) &= \frac{1}{2} x_N^\top Q_N x_N + q_N^\top x_N,
\end{align*}
where $q_k$, $r_k$, $Q_k$, $M_k$, $R_k$, and $q_N$, $Q_N$ are running and terminal cost weight matrices. 
Moreover, the stage-wise constraints~\eqref{eq:r_constraint} and~\eqref{eq:t_constraint} are characterized by matrices $D_{x,k} \in \mathbb{R}^{n_{c,k} \times n_x}$, $D_{u,k} \in \mathbb{R}^{n_{c,k} \times n_u}$, $e_k \in \mathbb{R}^{n_{c,k}}$, and a (Cartesian product of) convex cone(s) $\mathcal{K}_k$.
We solve \eqref{eq:conic_lqr} via an ADMM-based approach, where the core idea is to separately handle the dynamics and conic constraints to simplify the optimization process, as described below.

\subsection{Conic ADMM}
Following Garstka et al.~(\citeyear{garstka_cosmo_2021}), we reformulate~\eqref{eq:conic_lqr} as: 
\begin{subequations} \label{eq:admm_lqr}
\begin{align} 
\min_{\tilde{\boldsymbol{x}}, \tilde{\boldsymbol{u}}, \boldsymbol{x}, \boldsymbol{u}, \boldsymbol{s}}\; & \sum_{k=0}^{N-1}  \ell_k(\tilde{x}_k, \tilde{u}_k) + \ell_N(\tilde{x}_N) +  \sum_{k=0}^{N} I_{\mathcal{K}_k}(s_k),\\ \label{eq:cost_w_indicator}
\text{s.t.}\quad & \tilde{x}_{k+1} = A_k \tilde{x}_k + B_k \tilde{u}_k + c_k,\\
& D_{x,k} \tilde{x}_k + D_{u,k} \tilde{u}_k + s_k = e_k, \label{eq:stage_eq}\\
& \tilde{x}_k = x_k, \tilde{u}_k = u_k,\quad\quad k = [0,N-1], \label{eq:stage_consensus_eq} \\
& D_{x,N} \tilde{x}_N + s_N = e_N, \label{eq:terminal_eq} \\ 
& \tilde{x}_N = x_N, \label{eq:terminal_consensus_eq}
\end{align}
\end{subequations}
where $\tilde{x}_{k}$ and $\tilde{u}_{k}$ are auxiliary variables achieving consensus with $x_{k}$ and $u_{k}$, respectively, through \eqref{eq:stage_consensus_eq} and \eqref{eq:terminal_consensus_eq}. The indicator function $I_{\mathcal{K}}$ of the set $\mathcal{K}$ in~\eqref{eq:cost_w_indicator} integrates the conic constraints into the objective.
The ADMM algorithm solves problem~\eqref{eq:admm_lqr} through the following three-step iteration \citep{boyd2011distributed}:
\begin{subequations}
\begin{align}
    \boldsymbol{\tilde{x}}^{j+1}, \boldsymbol{\tilde{u}}^{j+1} &:= \argmin_{\boldsymbol{\tilde{x}}^, \boldsymbol{\tilde{u}}}~\mathcal{L}_{\rho,\sigma}(\boldsymbol{x}^j, \boldsymbol{u}^j, \boldsymbol{s}^j,\boldsymbol{y}^j,\boldsymbol{\tilde{x}},\boldsymbol{\tilde{u}}),\label{ADMM1}\\
    \boldsymbol{s}^{j+1} &:= \Pi_{\mathcal{K}} \left( \boldsymbol{\tilde{s}}^{j+1} + \rho^{-1} \boldsymbol{y}^j \right),\label{ADMM2}\\
    \boldsymbol{y}^{j+1}&:= \boldsymbol{y}^{j} + \rho(\boldsymbol{\tilde{s}}^{j+1}-\boldsymbol{s}^{j+1}), \label{ADMM3}
\end{align}
\label{eq:admm_steps}%
\end{subequations}
where \eqref{ADMM1} represents an unconstrained LQ problem that minimizes the augmented Lagrangian $\mathcal{L}_{\rho,\sigma}$ over $\tilde{\boldsymbol{x}}$ and $\tilde{\boldsymbol{u}}$, subject to linear dynamics.
The expression for the augmented Lagrangian is given by:
\begin{align} \label{eq:al_cost}
    \mathcal{L}_{\rho,\sigma}(\cdot) ~& = \sum_{k=0}^{N-1}  \ell_k(\tilde{x}_k, \tilde{u}_k) + \ell_N(\tilde{x}_N) +  \sum_{k=0}^{N} I_{\mathcal{K}_k}(s_k), \\
    ~&+ \dfrac{\rho}{2} \left\Vert D_{x,N} \tilde{x}_N + s_N - e_N + \rho^{-1}y_N \right\Vert_2^2 \nonumber \\
    ~&+ \sum_{k=0}^{N-1} \dfrac{\rho}{2} \left\Vert D_{x,k} \tilde{x}_k + D_{u,k} \tilde{u}_k + s_k - e_k + \rho^{-1}y_k  \right\Vert_2^2 \nonumber \\
    ~&+ \sum_{k=0}^{N} \frac{\sigma}{2} \left\Vert x_k - \tilde{x}_k \right\Vert_2^2 + \sum_{k=0}^{N-1} \frac{\sigma}{2} \left\Vert u_k - \tilde{u}_k \right\Vert_2^2, \nonumber
\end{align}
where the parameters $\rho > 0$ and $\sigma > 0$ act as regularization or step-size terms, and $\boldsymbol{y}$ is the vector of dual variables associated with the stage-wise constraints~\eqref{eq:stage_eq}.
\eqref{ADMM2} projects $\boldsymbol{s}$ into $\mathcal{K}$, and~\eqref{ADMM3} updates the dual variables. The $j$ superscript denotes the ADMM iteration index. 
An intermediate constraint relaxation step is often included between~\eqref{ADMM1} and ~\eqref{ADMM2}, which for the running variables takes the following form:
\begin{subequations}
\begin{align}
    (&{x_k}^{j+1}, {u_k}^{j+1}) := \alpha(\tilde{x}_k^{j+1}, \tilde{u}_k^{j+1}) + (1-\alpha)(x_k^j, u_k^j) \\
    &\tilde{s}_k^{j+1} := \alpha( e_k - D_{x,k} \tilde{x}_k^{j+1} - D_{u,k} \tilde{u}_k^{j+1}) + (1-\alpha) s_k^j,
\end{align}
\label{eq:admm_relaxation_steps}%
\end{subequations}
where $\alpha\in(0,2)$ is the relaxation parameter. 
If $\alpha>1$, \eqref{eq:admm_relaxation_steps} imposes over-relaxation, with demonstrated improved convergence using $\alpha$ values in the range $1.5-1.8$~\citep{eckstein1992douglas}.
In our implementation, the relaxation parameter is set to $\alpha = 1.6$, following ~\cite{osqp}.

\subsection{The LQ Problem}
In~\eqref{ADMM1}, we solve an LQ problem, which arises not only in ADMM-type methods but is also compatible with alternative constraint-handling approaches, including interior-point methods \citep{Domahidi_IPM_2012} and proximal augmented Lagrangian methods \citep{Lowenstein_alm_ocp_2024}.
To match the quadratic form of the cost function in~\eqref{eq:conic_lqr}, we reformulate the augmented Lagrangian in \eqref{eq:al_cost} as a function of terms of the following cost matrices:
\begin{subequations} \label{eq:cost_matrices}
\begin{align}
    \hat{Q}_k &= Q_k + \sigma I + \rho D_{x,k}^\top D_{x, k},  \qquad {\color{blue}\texttt{syrk}\; n_x^2 n_c}\\
    \hat{R}_k &= R_k + \sigma I + \rho D_{u,k}^\top D_{u, k},  \qquad {\color{blue}\texttt{syrk}\; n_u^2 n_c} \\
    \hat{M}_k &= M_k + \rho D_{u, k}^\top D_{x, k}, \;\;\qquad {\color{blue}\texttt{gemm}\; 2 n_x n_u n_c} \\
    \hat{q}_k &= q_k + D_{x, k}^\top \left(\rho (s_k^j - e_k) - y^j \right) - \sigma x^j,\\
    \hat{r}_k &= r_k + D_{u, k}^\top \left(\rho (s_k^j - e_k) - y^j \right) - \sigma u^j.
\end{align}
\end{subequations}
We highlight the linear algebra routines and floating-point operation (flop) counts in blue for reference, with additional context provided in the following section.

%%%%%%%%%%%%%%%%%%%%%%%%%%%%%%%%%%%%%%%%%%%%%%%%%%%%%%%%%%%
%%%%%%%%%%%%%%%%%%%%%%%%%%%%%%%%%%%%%%%%%%%%%%%%%%%%%%%%%%%
%%%%%%%%%%%%%%%%%%%%%%%%%%%%%%%%%%%%%%%%%%%%%%%%%%%%%%%%%%%
%%%%%%%%%%%%%%%%%%%%%%%%%%%%%%%%%%%%%%%%%%%%%%%%%%%%%%%%%%%

\begin{figure}[!t]
    \centering
    \includegraphics[width=0.9\linewidth]{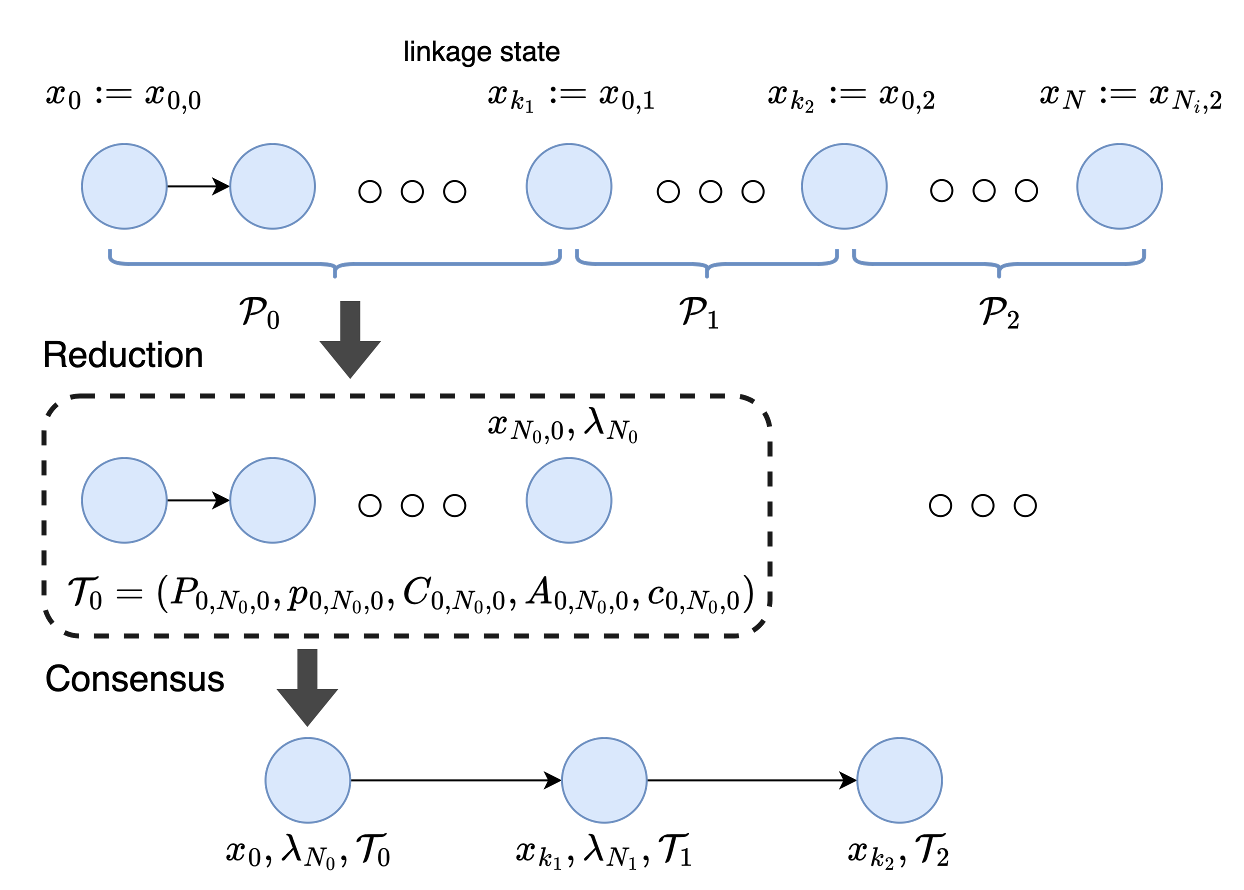}
    \caption{Parallel-in-time illustration for an LQ problem partitioned into $J=3$ subproblems. 
    Adjacent subproblems are coupled through a shared linkage state. Each subproblem in the reduction phase is a fixed-end LQ problem. The tuple $\mathcal{T}_i$, represents the components of the conditional value function at stage $i$. In the consensus phase, we combine the solutions of the subproblems and obtain the linkage states and Lagrange multipliers. Finally, the state and control input sequences are recovered.} 
    \label{fig:parlqr}
\end{figure}

\section{Parallel LQ Solver} \label{sec:parLQR}
Our parallel algorithm to address~\eqref{ADMM1} comprises a backward pass and a forward pass. The backward pass includes two phases, namely the reduction and consensus phases, as illustrated in Figure \ref{fig:parlqr}. 
The forward pass involves a parallel rollout of the linear dynamics.
This section focuses on the methods applied in the backward pass. 

We adopt a divide-and-conquer strategy to parallelize the computation of \eqref{ADMM1}.
In particular, we break down the LQ problem along the prediction horizon into $J$ similar subproblems $\{ \mathcal{P}_i \}_{i=0}^{J-1}$. 
Let $N_i$ and $k_i$ denote the horizon length and starting stage of the $i$-th subproblem.
We then define the state and input sequences for subproblem $i$ as:
\begin{subequations}
\begin{align}
    \boldsymbol{x}_i &= (x_{1,i}, \dots, x_{N_i, i}) := (x_{k_i+1}, \dots, x_{k_i + N_i}) \\
    \boldsymbol{u}_i &= (u_{0,i}, \dots, u_{N_i-1, i}) := (u_{k_i}, \dots, u_{k_i + N_i - 1}).
\end{align}
\end{subequations}
Two consecutive subproblems, $\mathcal{P}_i$ and $\mathcal{P}_{i+1}$, are linked by the condition $x_{N_i, i} = x_{0, i+1}$, ensuring continuity between the final state of $\mathcal{P}_i$ and the initial state of $\mathcal{P}_{i+1}$.
Therefore, in addition to the initial state $x_{0,i}$, we also parameterize each subproblem by its terminal state $x_{N_i, i}$, resulting in a fixed-end LQ problem $\mathcal{P}_i(x_{0,i}, x_{N_i, i})$:
% \begin{subequations} 
\begin{align} \label{eq:subproblem}
    \min_{\boldsymbol{\hat{x}}_i, \boldsymbol{u_i}} &\sum _{n=0}^{N_i-1} \ell_{n,i}(x_{n,i}, u_{n,i}) \\
    \text{s.t.} &\; x_{n+1,i} = A_{n, i} x_{n, i} + B_{n, i} u_{n, i} + c_{n, i},\, n \in [0, N_i-1], \nonumber %\\
\end{align} 
% \end{subequations}
where $\boldsymbol{\hat{x}}_i = \left( x_{1,i}, \dots, x_{N_i-1,i} \right)$ is the state sequence, excluding the initial and terminal states. 

In the remainder of this section, we describe the stages of our algorithm in detail and occasionally omit the subscript $i$ and the hat symbol for notational simplicity.

\subsection{Reduction Phase} \label{subsec:reduction_phase}
By dualizing the dynamics constraint coupling the terminal and second last stages, we reformulate \eqref{eq:subproblem} as follows:
\begin{align} \label{eq:fixed-end_LQR_dual} %
    \max_{\lambda_{N}}
    \min _{\boldsymbol{\hat{x}}, \boldsymbol{u}}\; & \sum_{n=0}^{N-1} \ell_{n} (x_{n}, u_{n}) \\
    & + \lambda_{N}^\top \left( A_{N-1}x_{N-1} + B_{N-1}u_{N-1} - x_{N} \right) \nonumber\\  
    \text{s.t.} &\; x_{n+1} = A_n x_n + B_n u_n + c_n, \; n \in [0, N-2] \nonumber,
\end{align} 
where $x_0$ and $x_N$ are treated as parameters. The inner minimization problem in \eqref{eq:fixed-end_LQR_dual} resembles a standard LQ problem with no terminal cost, but includes an additional term arising from the terminal state constraint.
We then apply dynamic programming to solve the inner optimization problem, starting from the minimization problem at stage $N-1$, where $g_{N-1, N}$ denotes the conditional dual value function from stage $N-1$ to $N$:
\begin{equation}
\begin{aligned}
    &g_{N-1, N} (x_{N-1}, x_{N}, \lambda_{N}) = \min_{u_{N-1}}\; \ell_{N-1} (x_{N-1}, u_{N-1}) \\
    &\quad\quad\quad\quad\quad +\lambda_{N}^\top \left( A_{N-1}x_{N-1} + B_{N-1}u_{N-1} - x_{N} \right).
\end{aligned}
\end{equation}
By analytically solving the unconstrained problem above, we derive the mathematical expression for the conditional dual value function from stage $n$ to $N$ as follows:
\begin{equation} 
\begin{aligned}
    &g_{n, N} (x_n, x_{N}, \lambda_{N}) = \beta + \frac{1}{2} x_{n}^\top P_{n, N} x_n + p_{n, N}^\top x_n \\
    &\quad\quad - \frac{1}{2} \lambda_{N}^\top C_{n, N} \lambda_{N} + \lambda_{N}^\top \left( A_{n, N}x_n + c_{n, N} - x_{N} \right),
\end{aligned}
\end{equation}
where $\beta$ is a constant value.
Given $g_{n+1, N}$, the Bellman equation at stage~$0 \leq n < N - 1$ is 
\begin{equation} \label{eq:stage_optimization}
    \min_{u_{n}}\; \ell_n(x_n, u_n) + g_{n+1, N}(A_n x_n + B_n u_n + c_n, x_N, \lambda_N). 
\end{equation}

We represent the action-value function in \eqref{eq:stage_optimization} by $Q_{n}$~\eqref{eq:Q_fun}, and derive the cost Hessians and gradients by grouping relevant terms as shown in~\eqref{eq:Q}:
\begin{equation} \label{eq:Q_fun}
    \begin{aligned}
        &Q_{n}(x_n, u_n, x_N, \lambda_N) = \beta + \begin{bmatrix}
            Q_{x, n} \\ Q_{u, n} \\ Q_{\lambda, n} 
        \end{bmatrix}^\top 
        \begin{bmatrix}
            x_n \\ u_n \\ \lambda_N 
        \end{bmatrix} \\
        &\qquad\qquad
        + \dfrac{1}{2}
        \begin{bmatrix}
            x_n \\ u_n \\ \lambda_N
        \end{bmatrix}^\top
        \begin{bmatrix}
            Q_{xx, n} & Q_{ux, n}^\top & Q_{\lambda x, n}^\top \\
            Q_{ux, n} & Q_{uu, n} & Q_{\lambda u, n}^\top \\
            Q_{\lambda x, n} & Q_{\lambda u, n} & Q_{\lambda\lambda, n} \\
        \end{bmatrix}
        \begin{bmatrix}
            x_n \\ u_n \\ \lambda_N 
        \end{bmatrix},
    \end{aligned}
\end{equation}
\begin{equation} \label{eq:Q}
\begin{aligned}
    Q_{xx, n} &= Q_n + A_n P_{n+1, N} A_n, \\
    Q_{uu, n} &= R_n + B_n P_{n+1, N} B_n, \\
    Q_{ux, n} &= M_n + B_n^\top P_{n+1, N} A_n, \\
    Q_{\lambda\lambda, n} &= -C_{n+1, N}, \\
    Q_{\lambda x, n} &= A_{n+1, N} A_n, \\
    Q_{\lambda u, n} &= A_{n+1, N} B_n, \\
    Q_{x, n} &= q_n + A_n^\top \left( P_{n+1, N} c_n + p_{n+1, N} \right), \\
    Q_{u, n} &= r_n + B_n^\top \left(P_{n+1, N} c_n + p_{n+1, N} \right), \\
    Q_{\lambda, n} &= A_{n+1, N} c_n + c_{n+1, N}.
\end{aligned}
\end{equation}

Next, by setting the gradient of $Q_{n}$ to zero, we solve the minimization problem~\eqref{eq:stage_optimization} to obtain the optimal control input $u_n^*$ and the conditional dual value function $g_{n, N}$. 
% described by $\left(P_{n, N}, p_{n, N}, C_{n, N}, A_{n, N}, b_{n, N} \right)$. 
The optimal control input is then:
\begin{align*}
    u_n^* &= -Q_{uu, n}^{-1} Q_{ux, n} x_n  -Q_{uu, n}^{-1} Q_{\lambda u, n}^\top \lambda_N - Q_{uu, n}^{-1} Q_{u, n} \\
    &:= K_n x_n + G_n \lambda_N + d_n.
\end{align*}
Substituting $u_n^*$ into the action-value function in \eqref{eq:Q_fun}, we obtain the key components of $g_{n, N}$:  
\begin{equation*}
\begin{aligned}
    P_{n,N} &= Q_{xx,n} - Q_{ux,n}^\top Q_{uu,n}^{-1} Q_{ux,n} \\
    C_{n, N} &= - Q_{\lambda\lambda,n} + Q_{\lambda u,n}\,Q_{uu,n}^{-1} Q_{\lambda u,n}^\top \\
    A_{n, N} &= Q_{\lambda x,n} - Q_{\lambda u,n}\,Q_{uu,n}^{-1} Q_{ux,n} \\
    p_{n,N} &= Q_{x,n} - Q_{ux,n}^\top Q_{uu,n}^{-1} Q_{u,n}  \\
    c_{n, N} &= Q_{\lambda,n} - Q_{\lambda u,n}\,Q_{uu,n}^{-1} Q_{u,n}. \\
\end{aligned}
\end{equation*}
Solving the first-stage minimization problem, we obtain the conditional value function $V_{i}$ \citep{SarkkaTemporal2023} for subproblem $i$ by maximizing the conditional dual value function over the dual variable $\lambda_{N_i}$:
\begin{equation} \label{eq:cond_val}
\begin{aligned}
    V_{i} (x_{0, i}, x_{N_i, i}) &= \max_{\lambda_{N_i}}\, g_{0, N_i, i} (x_{0, i}, x_{N_i, i}, \lambda_{N_i, i}),
\end{aligned} 
\end{equation}
which represents the optimal cost of the trajectory, conditioned on the initial and terminal states.
If the terminal state is not reachable, the value function tends toward infinity.
To simplify the notation, we introduce the following:
\begin{equation*}
\begin{aligned} 
    &\mathcal{T}_{i} := \left(P_{i,i+1}, p_{i,i+1}, C_{i,i+1}, A_{i,i+1}, b_{i,i+1}\right) := \\
    &\qquad\;\left(P_{0, N_i, i}, p_{0, N_i, i}, C_{0, N_i, i}, A_{0, N_i, i}, b_{0, N_i, i}\right),\; i \in [0, J-2].
\end{aligned}
\end{equation*}
For the last subproblem, we define $\mathcal{T}_{J-1} := \left(P_{J-1}, p_{J-1}\right)$.

The full reduction phase is shown in Algorithm~\ref{algo:reduction_sqrt}, with the linear algebra operation and the flop counts highlighted in blue, including only the cubic and quadratic terms. 
In Algorithm \ref{algo:reduction_sqrt} (Lines 7-11), instead of directly implementing \eqref{eq:Q}, we adopt the square-root Riccati recursion proposed in~\citep{Frison_LQ_2013}, where the recursion is expressed in terms of the Cholesky factor $L_{nn, x}$ rather than $P_n$. 
This method reduces the flop counts and improves spatial locality in the cache by packing the matrices.
The computational cost of the stage factorization is summarized in Table~\ref{tab:phase_cost}.
We highlight in bold the additional floating-point operations required by the stage factorization in the fixed-end reduction.
Moreover, during the reduction phase, we perform $J$ backward passes in parallel.
If the whole LQ problem is evenly divided, the computational load differs between the first $J-1$ backward passes and the final one. 
Specifically, the last backward pass employs the standard Riccati recursion, by excluding Lines 17-25.
To minimize unnecessary lag time, we balance the computational load by increasing the horizon length of the final subproblem.
Let $\psi$ denote the ratio between the computation time of the stage factorization in the fixed-end reduction (Lines 7-25) and that in the free-end reduction (Lines 7-16). 
Assuming that all CPUs have identical computing capabilities, we determine the horizon length $N^\prime$ of the first $J-1$ subproblems as follows, by equating the computation time of the first subproblem with that of the last one:
\begin{align} \label{eq:N_prime}
    N^\prime &\approx \texttt{int} \left( N / (\psi + J -1) \right),
\end{align}
which we denote as a load-balancing scheme.

\begin{table}[t]
\caption{\textsc{Cost of factorization for reduction and consensus phases}}
\label{tab:phase_cost}
\centering
\renewcommand{\arraystretch}{1.0}
% \small
\resizebox{0.95\linewidth}{!}{%
\begin{tabular}{@{}lll@{}}  % <-- now has three columns
\toprule
\textbf{Phase} & \textbf{Operation} & \textbf{Cost per stage (flops)} \\ \midrule

Reduction &
Factorization (free-end) &
$\begin{aligned}[t]
&\frac{7}{3} n_x^3 + 4n_x^2 n_u + 2n_x n_u^2 + \frac{1}{3} n_u^3 \\
&\quad + (n_x + n_u)^2 n_{c, k}
\end{aligned}$ \\

&
Factorization (fixed-end) &
$\begin{aligned}[t]
& \text{Free-end Factorization Cost} \\
&\quad + \mathbf{2n_x^3 + 6n_x^2 n_u + 2n_x n_u^2}
\end{aligned}$\\ \midrule

Consensus &
Factorization (LU) &
$\begin{aligned}[t]
26n_x^3/3
\end{aligned}$ \\

&
Factorization (Cholesky) &
$\begin{aligned}[t]
20n_x^3/3
\end{aligned}$ \\

\bottomrule
\end{tabular}%
}
% \vspace{-1pt}
\end{table}

\begin{algorithm}[!t]
% \setstretch{1.1}
\caption{Reduction phase}
\label{algo:reduction_sqrt}
\begin{algorithmic}[1]
\vspace{1mm}
\IF{\textit{not last subproblem}}
    \STATE $L_{xx, N} \leftarrow 0,\; p_N \leftarrow 0, C_N \leftarrow 0, F_N \leftarrow I$
\ELSE
    \STATE $P_N \leftarrow Q_n^{1/2},\; p_N \leftarrow q_N$
\ENDIF
\FOR{$n = N-1 \rightarrow 0$}  %\COMMENT{Backward Pass}
    \STATE $\mathcal{V}_n \leftarrow 
    \begin{bmatrix}
        B_n^\top \\ A_n^\top 
    \end{bmatrix} L_{xx, n+1}$ \hfill
    {\color{blue}$\texttt{trmm}\; n_x^2 (n_x + n_u)$}
    
    \STATE $
    \mathcal{M}_n \leftarrow 
    \begin{bmatrix}
        R_n & * \\
        S_n^\top & Q_n
    \end{bmatrix} + \mathcal{V}_n \mathcal{V}_n^\top$ \\ 
    $ \mathcal{M}_n := 
    \begin{bmatrix}
        R_n + B_n^\top P_{n+1} B_n \\
        S_n^\top + A_n^\top P_{n+1} B_n & Q_n + A_n^\top P_{n+1} A_n
    \end{bmatrix}
    $ \\ 
    \hfill {\color{blue}$\texttt{syrk}\; n_x (n_x + n_u)^2$}

    \STATE $
    \begin{bmatrix}
        L_{uu, n} & * \\
        L_{xu, n} & L_{xx, n} 
    \end{bmatrix} \leftarrow 
    \textsc{choFact} \left( \mathcal{M}_n \right)
    % $ \\ \hfill $\text{potrf}\; \frac{1}{3} (n_x + n_u)^3$
     $\hfill {\color{blue}$\frac{1}{3} (n_x + n_u)^3$}
    \STATE
    \STATE $P_{n+1}c_n  \leftarrow L_{xx, n+1}\left( L_{xx, n+1}^\top c_n \right) + p_{n+1}$ \hfill {\color{blue}$\texttt{trmv}\;2 n_x^2$}
    \STATE $
    \begin{bmatrix}
        l_n \\ p_n 
    \end{bmatrix} \leftarrow
    \begin{bmatrix}
        r_n \\ q_n 
    \end{bmatrix} +
    \begin{bmatrix}
        B_n^\top \\ A_n^\top 
    \end{bmatrix}
    P_{n+1}c_n
    $ \hfill {\color{blue}$\texttt{gemv}\; 2 n_x (n_x + n_u)$}
    \STATE $l_n \leftarrow L_{uu, n}^{-1} l_n$ \hfill {\color{blue}$\texttt{trsv}\;n_u^2$}
    \STATE $p_n \leftarrow p_n - L_{xu, n} l_n$ \hfill {\color{blue}$\texttt{gemv}\;2 n_x n_u$}
    \IF{\textit{not last subproblem}}
    \STATE $K_n \leftarrow L_{uu, n}^{-\top} L_{xu, n}^\top$ 
    \hfill {\color{blue}$\texttt{trsv}\;n_x n_u^2$}
    \STATE $d_n \leftarrow L_{uu, n}^{-\top} l_n$ \hfill {\color{blue}$\texttt{trsv}\;n_u^2$}
    \STATE $[FB\; FA] \leftarrow [F_{n+1} B_n\; F_{n+1} A_n]$\hfill {\color{blue}$\texttt{gemm}\;2n_x^2 n_u + 2 n_x^3$}
    % \STATE $G_n \leftarrow \textsc{choSolve}(L_{uu, n}, -(FB)^\top)$\hfill $2n_x n_u^2$
    \STATE $\Lambda_n \leftarrow - L_{uu, n}^{-1} (FB)^\top$ \hfill {\color{blue}$\texttt{trsv}\; n_x n_u^2$}
    % \STATE $C_n \leftarrow C_{n+1} - FB G_n$\hfill $2n_x^2 n_u$
    \STATE $C_n \leftarrow C_{n+1} + \Lambda_n^\top \Lambda_n$ \hfill {\color{blue}$\texttt{gemm}\; 2n_x^2 n_u$}
    \STATE $F_n \leftarrow FA + FB K_n$\hfill {\color{blue}$\texttt{gemm}\; 2n_x^2 n_u$}
    \STATE $f_n \leftarrow f_{n+1} + F_{n+1} (B_n d_n + c_n)$ \hfill {\color{blue}$\texttt{gemv}\;2 n_x n_u + 2 n_x^2$}
    \ENDIF
\ENDFOR
\STATE $P_{0, N} \leftarrow L_{xx, 0} L_{xx, 0}^\top$, $p_{0, N} \leftarrow p_0$
\STATE $C_{0, N} \leftarrow C_0$, $A_{0, N} \leftarrow F_0$, $c_{0, N} \leftarrow f_0$
\end{algorithmic}
\end{algorithm}

\begin{algorithm}[!t]
% \setstretch{1.1}
\caption{Consensus phase}
\label{algo:consensus}
\begin{algorithmic}[1]
\vspace{1mm}
\STATE $P_{J-1} \leftarrow P_{N_{J-1}}, p_{J-1} \leftarrow p_{N_{J-1}}$, $x_{0}^{\text{lk}} \leftarrow x_0$ 
\FOR{$i = J-2 \rightarrow 0$}
    \STATE $CP \leftarrow C_{i,i+1}P_{i+1},\; PA \leftarrow P_{i+1}A_{i,i+1}$ \hfill {\color{blue}$\texttt{gemm}\; 4n_x^3$}
    \STATE $\text{lufact}_i \leftarrow \textsc{luFact}(I + CP)$ \hfill {\color{blue}$\frac{2}{3} n_x^3$} 
    \STATE $\mathcal{D}_i \leftarrow \textsc{luSolve}(\text{lufact}_i, A_{i,i+1})$ \hfill {\color{blue}$2n_x^3$}
    \STATE $P_i \leftarrow P_{i,i+1} + \mathcal{D}_i^\top P A$ \hfill {\color{blue}$\texttt{gemm}\;2n_x^3$}
    % \STATE $P_i = P_{i,i+1} + A_{i,i+1}^\top (I + P_{i+1}C_{i,i+1})^{-1} P_{i+1}A_{i,i+1}$
    \STATE $p_i \leftarrow p_{i,i+1} + \mathcal{D}_i^\top (P_{i+1} c_{i,i+1} + p_{i+1})$ \hfill {\color{blue}$\texttt{gemv}\; 4 n_x^2$}
\ENDFOR
\FOR{$i = 0 \rightarrow J-2$}
    \STATE $e_i \leftarrow A_{i,i+1} x_{i}^{\text{lk}} + c_{i,i+1} - C_{i,i+1} p_{i+1}$ \hfill {\color{blue}$\texttt{gemv}\; 4 n_x^2$}
    \STATE $x_{i+1}^{\text{lk}} \leftarrow \textsc{luSolve}(\text{lufact}_i, e_i) $ \hfill {\color{blue}$2 n_x^2$}
    \STATE $\lambda_{i+1}^{\text{lk}} \leftarrow P_{i+1} x_{i+1}^{\text{lk}} + p_{i+1}$ \hfill {\color{blue}$\texttt{gemv}\; 2 n_x^2$}
\ENDFOR
\end{algorithmic}
\end{algorithm}

\subsection{Consensus Phase}
Next, we aim to combine the subproblems by considering a minimization problem associated with all linkage states.
Let the linkage state sequence be denoted as:
$
\boldsymbol{x}^{\text{lk}} = (x^{\text{lk}}_1, \dots, x^{\text{lk}}_{J-1}) := (x_{k_1}, \dots, x_{k_{J-1}}).
$
The optimization problem is then formulated as:  
\begin{equation} \label{eq:comb}
\begin{aligned}
    V(x_0) = \min_{\boldsymbol{x}^{\text{lk}}}\; V_{J-1}(x^{\text{lk}}_{J-1}) + \sum_{i=0}^{J-2} V_i(x^{\text{lk}}_{i}, x^{\text{lk}}_{i+1}),
\end{aligned}
\end{equation}
where $x^{\text{lk}}_0$ is defined as $x_0$.
By substituting \eqref{eq:cond_val} into \eqref{eq:comb}, we derive a min-max problem over the linkage states and corresponding dual variables: 
\begin{equation} \label{eq:comb_min_max} 
\begin{aligned}
    \min_{\boldsymbol{x}^{\text{lk}}}\,
    \max_{\boldsymbol{\lambda}^{\text{lk}}}\; V_{J-1}(x^{\text{lk}}_{J-1}) + \sum_{i=0}^{J-2} g_{0, N_i, i} (x^{\text{lk}}_{i}, x^{\text{lk}}_{i+1}, \lambda^{\text{lk}}_{i+1}),
\end{aligned}
\end{equation}
where $\boldsymbol{\lambda}^{\text{lk}} = (\lambda^{\text{lk}}_{1}, \dots, \lambda^{\text{lk}}_{J-1}) := (\lambda_{N_0}, \dots, \lambda_{N_{J-2}})$. 
The optimality conditions for \eqref{eq:comb_min_max} are obtained by taking the gradient of the objective function with respect to $\boldsymbol{x}^{\text{lk}}$ and $\boldsymbol{\lambda}^{\text{lk}}$, yielding the following system of equations:
\begin{align}
    -\lambda_{i-1}^{\text{lk}} + P_{i-1, i} x_{i-1}^{\text{lk}} + A_{i-1, i}^\top \lambda_{i}^{\text{lk}} + p_{i-1, i} &= 0,\; i \in [2, J-1] \nonumber\\
    A_{i-1,i} x_{i-1}^{\text{lk}} - C_{i-1,i} \lambda_i^{\text{lk}} - x_i^{\text{lk}} + c_{i-1, i} &= 0,\; i \in [1, J-1] \nonumber\\
    -\lambda_{J-1}^{\text{lk}} + P_{J-1} x^{\text{lk}}_{J-1} + p_{J-1} &= 0.
\end{align}
Leveraging the special structure of the resulting banded KKT matrix, we perform a backward elimination to establish the relationship between $x_i$ and $\lambda_i$:
\begin{align} \label{eq:consensus_backward}
    \lambda_i^{\text{lk}} &= P_i x_i^{\text{lk}} + p_i \\
    P_i &= P_{i,i+1} + A_{i,i+1}^\top (I + P_{i+1}C_{i,i+1})^{-1} P_{i+1}A_{i,i+1} \nonumber \\
    p_i &= p_{i,i+1} + A_{i,i+1}^\top (I + P_{i+1}C_{i,i+1})^{-1} (P_{i+1} c_{i,i+1} + p_{i+1}), \nonumber
\end{align}
where $P_i$ and $p_i$ are the parameters associated with the value function of an LQ problem. 
During the backward pass, the connection between two consecutive states is:
\begin{equation} \label{eq:consensus_forward}
    x_{i+1}^{\text{lk}} = (I + C_{i,i+1}P_{i+1})^{-1} (A_{i,i+1} x_{i}^{\text{lk}} + c_{i,i+1} - C_{i,i+1} p_{i+1}).
\end{equation}
Once $\{P_i\}$ and $\{p_i\}$ are computed, we use \eqref{eq:consensus_backward} and \eqref{eq:consensus_forward} to obtain the linkage states and Lagrangian dual variables.
The consensus phase is described in Algorithm~\ref{algo:consensus}.
Unlike \citep{jallet_parallel_2024}, our method only requires $P_{i,i+1}$ to be positive semi-definite rather than positive definite, at the cost of approximately $2n_x$ additional flop counts per stage.

\begin{figure*}[!t]
    \centering
    % --- Subfigure (a) ---
    \begin{subfigure}[t]{0.45\linewidth}
        \centering
        \includegraphics[width=\linewidth]{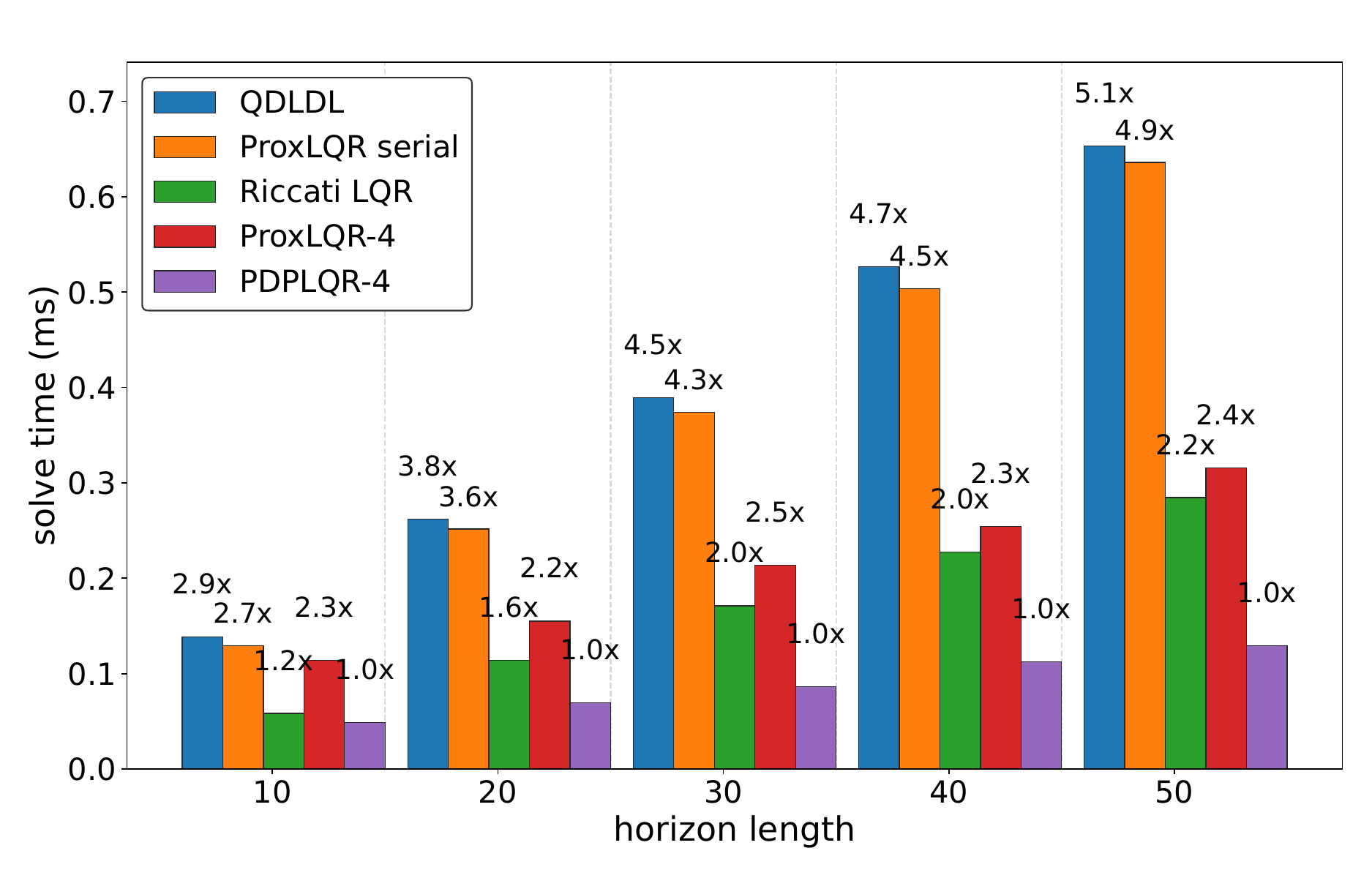}
        \vspace{-20pt}
        \caption{quadrotor hovering control}
        \label{fig:quadrotor_benchmark}
    \end{subfigure}
    \hfill
    % --- Subfigure (b) ---
    \begin{subfigure}[t]{0.45\linewidth}
        \centering
        \includegraphics[width=\linewidth]{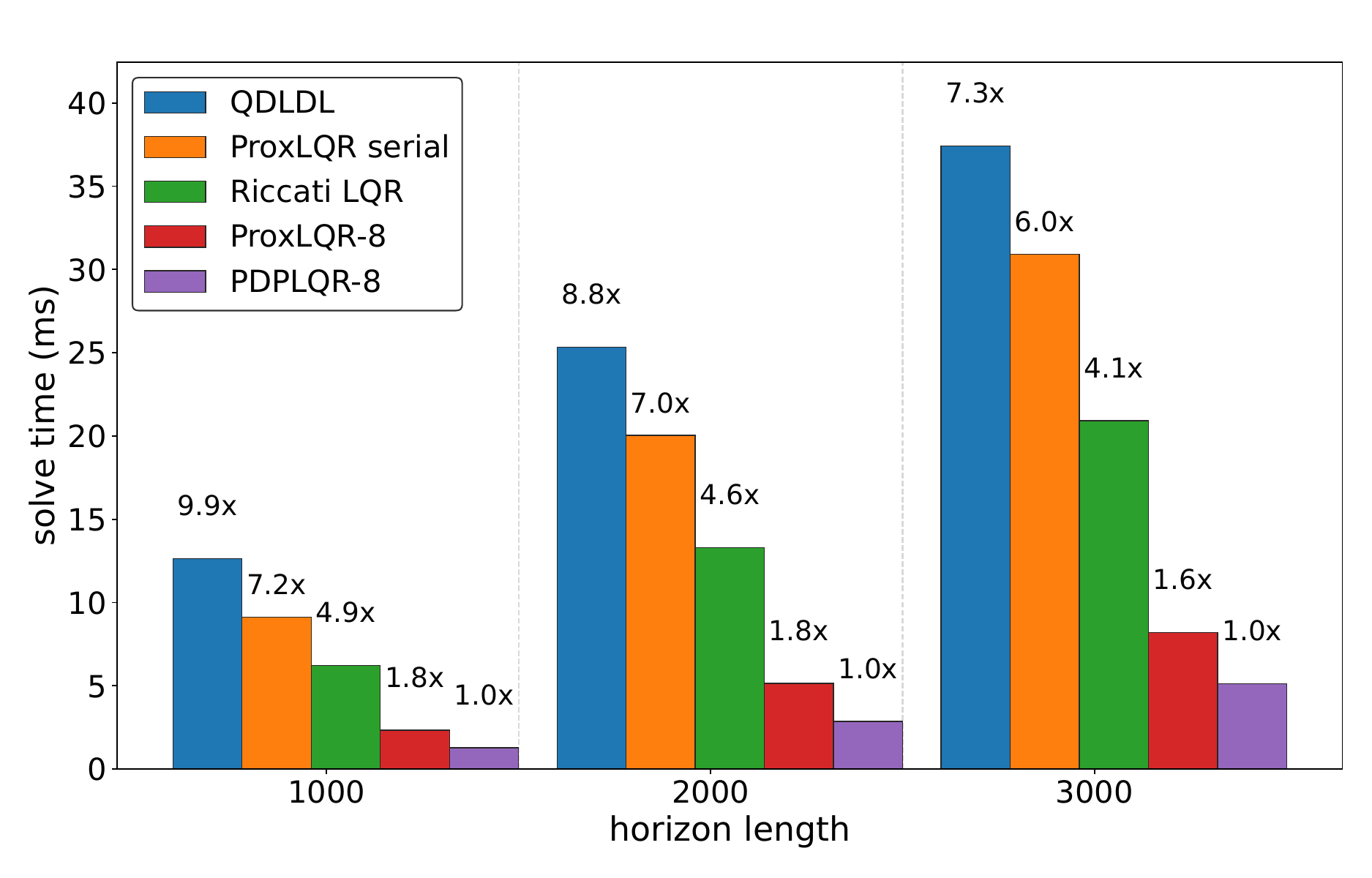}
        \vspace{-20pt}
        \caption{low-thrust orbital transfer planning}
        \label{fig:ksorbital_benchmark}
    \end{subfigure}
    % \vspace{-8pt}
    \caption{Benchmarking results across various horizon lengths. Above-bar numbers indicate the speedups achieved by our \texttt{PDPLQR} solver over realted work. In the quadrotor hovering control experiment, load balancing is turned off.}
    \label{fig:benchmark}
    % \vspace{-5mm}
\end{figure*}

\section{Numerical Results} \label{sec:results}
% \subsection{Methodology} 
The objectives of this numerical study are threefold: (1) to evaluate the computational performance of the proposed solver across a wide range of horizon lengths, (2) to examine the effects of varying the number of subproblems and stage constraints, and (3) to validate the effectiveness of the proposed load-balancing scheme via \eqref{eq:N_prime}.
In addition, we compare the performance of our parallel solver \texttt{PDPLQR} against several high-performance, state-of-the-art, LQ solvers: one based on the square-root Riccati recursion \citep{Frison_LQ_2013} adopted in acados \citep{Verschueren2021} and identified as \texttt{Riccati LQR}, the proximal LQ solver (\texttt{ProxLQR}) integrated in Aligator \citep{jallet_parallel_2024}, and
the widely used sparse linear solver \texttt{QDLDL} \citep{osqp}.

All simulations are conducted, using Google Benchmark, on a laptop with a \SI{2.30}{\giga\hertz} Intel Core i7-11800H processor (8 physical cores) and \SI{16}{\giga\byte} RAM. 
We implement our parallel method in C++ using \texttt{Eigen} for dense linear algebra. 
% We employ \texttt{Eigen}'s $\texttt{triangularView}$ for improved computational efficiency. 
The OpenMP API is used to enable parallel execution. 
To minimize the overhead of online thread creation, all threads are created and bound to specific CPU cores during the problem setup.
We disable Turbo Boost for stable thermal performance during benchmarking.
Moreover, for a fair comparison, we use Aligator v0.16.0, which is optimized for explicit dynamics.

\subsection{Benchmarking Results}
For benchmarking, we focus on two real-world engineering problems: (a) quadrotor hovering control representing a short-horizon case and (b) low-thrust orbital transfer planning representing a long-horizon case. 
The quadrotor model is linearized around the hovering condition, with 12 states $(n_x = 12)$ and 4 control inputs $(n_u = 4)$. 
We enforce box constraints on the state and input variables, resulting in $n_{c,k} = 16$ constraints per stage.
The orbital dynamics model is adopted from \citep{TracyManchester2021_LowThrustKS}, consisting of 13 states and 2 control inputs. The imposed stage constraint requires that the $\ell$-2 norm of the control input vector $u_k$ remains below a specific bound, $\Vert u_k \Vert_2 \leq u_{\max}$. This is a typical second-order cone constraint.

Figure \ref{fig:benchmark} illustrates the impact of various horizon lengths. \texttt{PDPLQR-J} refers to the proposed solver with $J$ subproblems. 
The solve time includes the backward and forward passes required for solving one LQ problem (the first ADMM step).
The proposed solver, \texttt{PDPLQR}, demonstrates superior performance on the quadrotor hovering control problem across all tested horizon lengths, achieving a speedup exceeding $4\times$ over \texttt{QDLDL} and $2\times$ over \texttt{ProxLQR} for $N > 20$. 
\texttt{ProxLQR} performs an $LDL^\top$ factorization of a matrix of size $(n_u + n_{c,k}) \times (n_u + n_{c,k})$, instead of a Cholesky factorization of a matrix of size $n_u \times n_u$. This approach is more numerically stable, but comes at the cost of increased computation time. 
In the low-thrust planning problem, we partition the LQ problem into eight subproblems to achieve better performance.
The parallel solvers, \texttt{PDPLQR} and \texttt{ProxLQR}, consistently outperform the serial implementations; moreover, our solver is more than $35\%$ faster than the \texttt{ProxLQR}.

\begin{figure}[!t]
    \centering
    \includegraphics[width=0.90\linewidth]{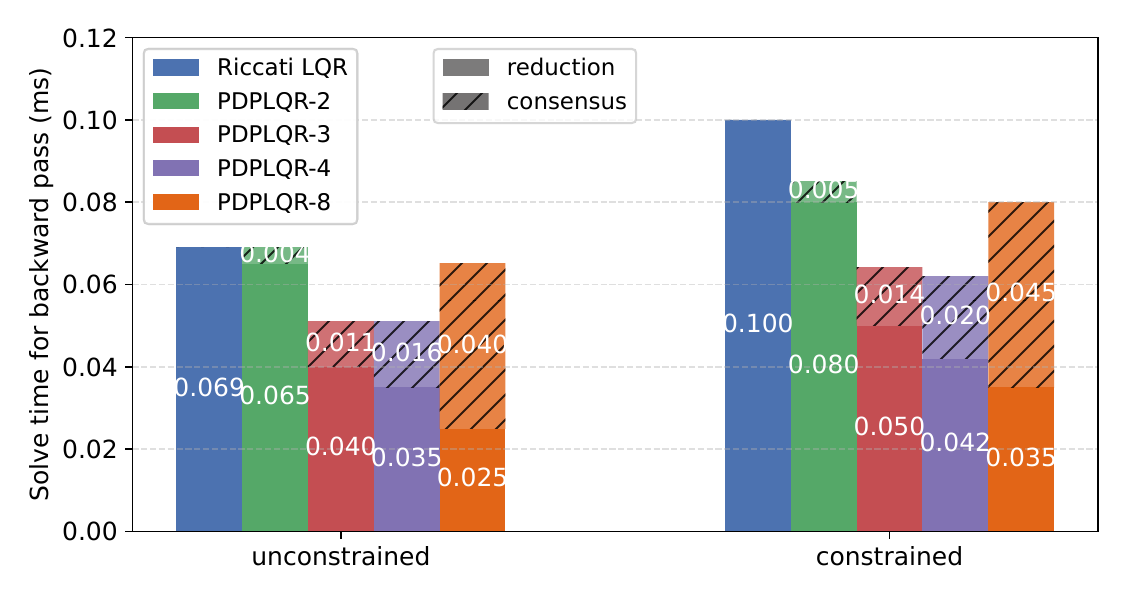}
    \caption{Solve times for the backward pass on the quadrotor control problem under unconstrained and constrained settings. Inside-bar numbers indicate phase solve times.}
    \label{fig:uncon_con}
\end{figure}

\subsection{Ablation Studies}
We now delve into assessing the backward pass (Algorithm \ref{algo:reduction_sqrt}), which is the most computationally demanding component of the LQ-problem solution process.  
Table \ref{tab:phase_cost} shows that solving the fixed-end LQ problem is more computationally expensive than its free-end counterpart.
By substituting $n_x= 12, n_u = 4, n_{c} =0$ into Table \ref{tab:phase_cost}, we can estimate that the ratio between the flop counts of the fixed-end and free-end stage factorizations is $2.08$. 
This value implies that splitting the horizon evenly into two parts offers little computational advantage. 
To validate this finding, we measure the computation time of the backward pass for the \texttt{Riccati LQR} and our parallel solver, \texttt{PDPLQR}. The results for the problem with $N = 20$ are shown in Figure \ref{fig:uncon_con}. 
As expected, dividing the LQ problem into three subproblems reduces the computation time. 
However, we note that it is not beneficial to use an excessive number of threads for short-horizon problems $(N \leq 20)$ since the factorization in the consensus phase often incurs a higher computational cost than that in the reduction phase, as indicated in Table \ref{tab:phase_cost}.
Next, we turn to the constrained LQ problem. 
As illustrated on the right side of Figure \ref{fig:uncon_con}, all parallel solvers outperform the \texttt{Riccati LQR}, since the flops associated with the stage constraints, $(n_x + n_u)^2 n_{c}$, offset the additional computational cost introduced by the fixed-end formulation. 
Overall, our method is effective for short-horizon problems with a large number of constraints and when properly choosing the number of segments.

Next, we investigate the effect of the load-balancing technique on the long-horizon problem introduced in Section \ref{subsec:reduction_phase}.
With the horizon length fixed at $2000$, we vary the number of subproblems.
Without load balancing, the solve times (\textrm{ms}) were 11.26, 5.61, 3.88, and 3.01 for $J = 2$, $4$, $6$, and $8$, respectively. Incorporating load balancing reduced these to 8.81, 4.95, 3.56, and \SI{2.87}{ms}, corresponding to improvements of $21.7\%$, $11.7\%$, $8.1\%$, and $4.9\%$.
The load-balancing effect becomes negligible, as the value of $N^\prime$~\eqref{eq:N_prime} is close to the evenly divided horizon length.

\section{Conclusion} \label{sec:conclusion}
In this work, we demonstrate that temporal parallelization on multi-core CPUs enables efficient solutions to LQ problems for optimal control problems with various horizon lengths.
In the future, we plan to extend our open-source solver to support additional constraint-handling methods, including proximal augmented Lagrangian methods and interior-point methods, and to adapt it to other parallel computing architectures like GPUs.

\section{DECLARATION OF GENERATIVE AI AND AI-ASSISTED TECHNOLOGIES IN THE WRITING PROCESS}
ChatGPT helped to enhance the writing quality of this work. Nonetheless, the authors reviewed and edited the manuscript throughout the full writing process, and assume full responsibility for the content of this publication.

\small
\bibliographystyle{ifacconf}
\bibliography{ifacconf}

\end{document}